# STABILITY ANALYSIS OF QUATERNION-VALUED NEURAL NETWORKS with LEAKAGE DELAY AND ADDITIVE TIME-VARYING DELAYS

Qun Huang and Jinde Cao

School of Mathematics, Southeast University, Nanjing, China

*ABSTRACT*

*In this paper, the stability analysis of quaternion-valued neural networks (QVNNs) with both leakage delay and additive time-varying delays is proposed. By employing the Lyapunov-Krasovskii functional method and fully considering the relationship between time-varying delays and upper bounds of delays, some sufficient criteria are derived based on reciprocally convex method and several inequality techniques. The stability criteria are established in two forms: quaternion-valued linear matrix inequalities (QVLMIs) and complex-valued linear matrix inequalities (CVLMIs),in which CVLMIs can be directly resolved by the Yalmip toolbox in MATLAB. Finally, an illustrative example is presented to demonstrate the validity of the theoretical results.*

*KEYWORDS*

*Quaternion-valued Neural Networks, Stability Analysis, Lyapunov-KrasovskiiFunctional, Leakage Delay, Additive Time-varying Delays*

## 1. INTRODUCTION

In the past decades, real-valued neural networks (RVNNs) have been successfully appliedin secure communication, information processing, engineer optimization, automatic control engineering and other areas. Correspondingly, numerous meaningful results have been reported [1-6]. However, RVNNs have its own limitations, such as the detection of symmetry problem cannot be resolved by a real-valued neuron, while it can be well solved by a complex-valued neuron [7]. In addition, the problem involving with ultrasonic wave, electromagnetic processing, quantum wave can be also well resolved by the complex number. Therefore, the performance of complex-valued neural networks (CVNNs) is more preferable than that of RVNN sin practical application with complex signals, and CVNNs have captured plenty of attentions from different areas [8-9]. In the past few years, it has drawn considerable attention to the dynamics of complex-valued neural networks and there have been lots of significant results associated with such kind of topics, see [10-11] and the references cited therein.

The quaternions are members of a noncommutative division algebra invented independently by William Rowan Hamilton in 1843. Some operation laws such as the commutativtiy of multiplication are not yet applicable for quaternions, which is quite different from the real or complex numbers. Owing to this difficulty, the research of quaternion had almost remained stagnant for a long period of time in the past. Recently, the resurgence of the study for quaternion systems is underway and an increasing spectrum of applications based on quaternions are found in various fields, such as quantum mechanics, attitude control, computergraphics and signal processing [12-13]. It has been proven that neural networks along with quaternion possess better performances and wider applications than both RVNNs and CVNNs. Actually, the three-dimensional and four-dimensional data can be expressed as an entirety, which is more authentic and reliable in modeling of practical application, and quaternion-valued neural networks (QVNNs) emerge at the right moment. Nowadays, increasing scholars are dedicated to investigating the dynamical behaviors of QVNNs. For instance, some sufficient criteria were proposed in the form of LMIs to guarantee the µ-stability of QVNNs with unbounded and non-differentiable time-varying delays in [14] and [15], respectively. In [16], by employing matrix measure and Halanay inequality technique, the problem of global exponential stability for delayed QVNNs was addressed successfully. Chen and Song [17] concentrated on the robust stability issue for delayed QVNNs based on homeomorphism mapping theorem and inequality techniques. Furthermore, the stability issue for both continuous-time and discrete-time QVNNs was investigated in [18]. Besides, some algebraic conditions were established to guarantee the global dissipativity for delayed QVNNs [19].

Time delays are inevitable in neural system owing to the limited propagation velocity between different neurons. Dynamical behaviors of neural networks could become more complicated owing to the existence of time delays, and it may result in performance degradation, such as instability, oscillation, bifurcation and so forth. Usually, the time delay in the state is supposed to appear in a singular form. Nevertheless, Zhao et al. [20] demonstrated that signals transmissions may experience a few segments of networks in several practical situations and different conditions of network transmission probably result in successive delays with different properties. By applying the convex polyhedron method, a less conservative delay-dependent stability criterion was proposed in [21]. Tian and Zhong [22] conducted further investigation on this issue by constructing augmented Lyapunov-Krasovskii functional and employing the reciprocally convex method, which is initially proposed by Park et al. [23]. Liang et al. generalized the reciprocally convex method to the complex domain and investigated the state estimation problem for complex-valued neural networks with two additive time-varying delays [24]. To the best ofthe authors' knowledge, up to now, few scholars have taken the stability problem of quaternion-valued neural networks with additive time-varying delays into consideration.

Enlightened by the aforementioned discussions, the aim of this paper is to conduct the stability analysis for quaternion-valued neural networks with both leakage delay and two additive time-varying delays. The remainder of this paper is organized as follows. In Section 2, the model description, several necessary hypotheses and lemmas are given. In Section 3, some sufficient criteria for the global asymptotical stability of QVNNs are derived based on reciprocally convex method and several inequality techniques. In Section 4, an illustrative example is presented to validate the effectiveness of the obtained results. Finally, conclusions are drawn in Section 5.

Notations: Let $\mathbb{R}$, $\mathbb{C}$ and $\mathbb{Q}$ stand for the real field, the complex field and the skew field of quaternions, respectively. Let $\mathbb{R}^{m\times n}$, $\mathbb{C}^{m\times n}$ and $\mathbb{Q}^{m\times n}$ separately denote $m\times n$ matrices with entries from $\mathbb{R}$, $\mathbb{C}$ and $\mathbb{Q}$. The notations $A^T$, $\bar{A}$ and $A^*$ represent the transpose, the conjugate and the conjugate transpose matrix of $A$, respectively. $A$ is referred to as Hermitian if $A=A^*$. The notation $X\geq Y(X>Y)$ means that $X-Y$ is positive semidefinite (positive definite, respectively). Moreover, the notation $*$ denotes the conjugate transpose of an appropriate block in a Hermitian matrix, while the notation ■ denotes the negative transpose of an appropriate block in a skew-symmetric matrix.

## 2. PRELIMINARIES

### 2.1. Quaternion Algebra

The quaternion is an extension of the complex number, and a quaternion $m \in \mathbb{Q}$ can be described in the following form:

$$m = m_0 + m_1 i + m_2 j + m_3 k,$$

where $m_0, m_1, m_2, m_3 \in \mathbb{R}$. The quaternion imaginary units $i, j, k$ obey the following rules:

$$i^2 = j^2 = k^2 = ijk = -1,$$

$$ij = -ji = k, ik = -ki = j, jk = -kj = i.$$

From which one can note that the quaternion multiplication is not commutative.
We proceed to introduce some basic operations of quaternion algebra. For two quaternions $m = m_0 + m_1 i + m_2 j + m_3 k$ and $n = n_0 + n_1 i + n_2 j + n_3 k$, the sum and product of $m$ and $n$ are defined as:

$$m + n = (m_0 + n_0) + (m_1 + n_1)i + (m_2 + n_2)j + (m_3 + n_3)k,$$

And

$$mn = (m_0 n_0 - m_1 n_1 - m_2 n_2 - m_3 n_3) + (m_0 n_1 + m_1 n_0 + m_2 n_3 - m_3 n_2) \\ + (m_0 n_2 + m_2 n_0 - m_1 n_3 + m_3 n_1)j + (m_0 n_3 + m_3 n_0 + m_1 n_2 - m_2 n_1)k.$$

In addition, the conjugate transpose of $m$ is defined as $m^* = m_0 - m_1 i - m_2 j - m_3 k$. The modulus of $m$ is denoted by $|m|$ and denoted as

$$|m| = \sqrt{mm^*} = \sqrt{(m_0)^2 + (m_1)^2 + (m_2)^2 + (m_3)^2}.$$

### 2.2. Model Formulation and Basic Lemmas

Consider the following quaternion-valued neural networks with both leakage delay and additive time-varying delays:

$$\dot{y}(t) = -Cy(t - \delta) + Ag(y(t)) + Bg\left(y(t - d_1(t) - d_2(t))\right) + h(t), \quad (1)$$

where $y(t) = (y_1(t), y_2(t), \ldots, y_n(t))^T \in \mathbb{Q}^n$ denotes the state vector, $C = diag\{c_1, c_2, \ldots, c_n\} \in \mathbb{R}^{n \times n}$ with $c_i > 0$ is the self-feedback connection weight matrix for $i \in \{1, 2, \ldots, n\}$. $A, B \in \mathbb{Q}^{n \times n}$ are the interconnection matrices which stand for the weight coefficients of the neurons. $g(x(t)) = (g_1(y_1(t)), g_2(y_2(t)), \ldots, g_n(y_n(t)))^T \in \mathbb{Q}^n$ represents the neuron activation function at time $t$; $h(t) \in \mathbb{Q}^n$ denotes the external input vector; $\delta$ is referred to as the leakage delay which satisfies $\delta \geq 0$; $d_1(t)$ and $d_2(t)$ represent the two delay components in the state.

In order to simplify the model, we assume that $y^*$ is an equilibrium point for (1). By applying the transformation $x(t) = y(t) - y^*$, system (1) is further converted to:

$$\dot{x}(t) = -Cx(t-\delta) + Af(x(t)) + Bf(x(t-d_1(t)-d_2(t))), \quad (2)$$

where

$$x(t) = (x_1(t), x_2(t), \ldots, x_n(t))^T, f(x(\cdot)) = (f_1(x_1(\cdot)), f_2(x_2(\cdot)), \ldots, f_n(x_n(\cdot)))^T, f_i(x_i(\cdot)) = g_i(x_i(\cdot) + y_i^*) - g_i(y_i^*) \ (i = 1,2,\ldots,n).$$

Throughout this paper, several assumptions play a crucial role and thus are presents as:

*Assumption 1*: $d_1(t)$ and $d_2(t)$ are continuous functions and satisfy

$$0 \leq d_1(t) \leq d_1, \quad 0 \leq d_2(t) \leq d_2,$$
$$\dot{d}_1(t) \leq \mu_1, \quad \dot{d}_2(t) \leq \mu_2,$$

where $d_i$ and $\mu_i$ ($i = 1,2$) are positive constants.

*Assumption 2*: For any $j \in \{1,2,\ldots,n\}$, there exists a positive constant $\gamma_j$ such that

$$|f_j(u_1) - f_j(u_2)| \leq \gamma_j |u_1 - u_2|$$

for all $u_1, u_2 \in \mathbb{Q}$.

Subsequently, we denote $d(t) = d_1(t) + d_2(t), d = d_1 + d_2, \mu = \mu_1 + \mu_2$ and $\Gamma = diag\{\gamma_1, \gamma_2, \ldots, \gamma_n\}$ for simplicity. Before deriving the main results, the following lemmas are instrumental.

**Lemma 2.1.** Let $u, v \in \mathbb{Q}, A, B \in \mathbb{Q}^{n \times n}, C \in \mathbb{C}^{n \times n}$. Then

(i) $|u + v| \leq |u| + |v|$, and $|uv| \leq |u||v|$;
(ii) $(AB)^* = B^*A^*$;
(iii) $(AB)^{-1} = B^{-1}A^{-1}$, if $A, B$ are invertible;
(iv) $(A^*)^{-1} = (A^{-1})^*$, if $A$ is invertible;
(v) $u$ can be uniquely expressed as $u = u_1 + u_2 j$, where $u_1, u_2 \in \mathbb{C}$;
(vi) $jC = \bar{C}j$ and $jCj = -\bar{C}$.

**Remark 2.1.** According to property (5) in Lemma 2.1, any quaternion matrix $A \in \mathbb{Q}^{n \times n}$ can be uniquely expressed as $A = A_1 + A_2 j$, where $A_1, A_2 \in \mathbb{C}^{n \times n}$.

**Lemma 2.2.** Let $A = A_1 + A_2 j, B = B_1 + B_2 j$, where $A, B \in \mathbb{Q}^{n \times n}$ and $A_1, A_2, B_1, B_2 \in \mathbb{C}^{n \times n}$. Then

(i) $A^* = A_1^* - A_2^T j$;
(ii) $AB = (A_1 B_1 - A_2 \bar{B}_2) + (A_1 B_2 + A_2 \bar{B}_1) j$.

**Lemma 2.3.** Let $A \in \mathbb{Q}^{n \times n}$ be a Hermitian matrix and $A = A_1 + A_2 j$, where $A_1, A_2 \in \mathbb{C}^{n \times n}$. Then $A < 0$ is equivalent to

$$\begin{pmatrix} A_1 & -A_2 \\ \bar{A}_2 & \bar{A}_1 \end{pmatrix} < 0.$$

**Remark 2.2.** Lemma 2.3 reveals the equivalence between the negative definiteness of a $n \times n$ quaternion matrix and the negative definiteness of a $2n \times 2n$ complex matrix.

**Lemma 2.4**. Suppose $M \in \mathbb{Q}^{n \times n}$ is a positive definite Hermitian matrix and $\omega(s):[a,b] \to \mathbb{Q}^n$ is a vector valued function. If the integrations concerned are well-defined, then

$$\left(\int_a^b \omega(s)ds\right)^* M \left(\int_a^b \omega(s)ds\right) \leq (b-a)\int_a^b \omega^*(s)M\omega(s)ds.$$

**Lemma 2.5**. For any vector $\xi \in \mathbb{Q}^m$, $\alpha \in (0,1)$, a positive Hermitian matrix $P \in \mathbb{Q}^{n \times n}$, and matrices $W_1, W_2 \in \mathbb{Q}^{n \times m}$, define the function $\Xi(\alpha, P)$ as

$$\Xi(\alpha, P) = \frac{1}{\alpha}\xi^* W_1^* P W_1 \xi + \frac{1}{1-\alpha}\xi^* W_2^* P W_2 \xi.$$

If there exists a matrix $X \in \mathbb{Q}^{n \times n}$ satisfying

$$\begin{pmatrix} P & X \\ * & P \end{pmatrix} \geq 0,$$

then

$$\min_{\alpha \in (0,1)} \Xi(\alpha, P) \geq \begin{pmatrix} W_1 \xi \\ W_2 \xi \end{pmatrix}^* \begin{pmatrix} P & X \\ * & P \end{pmatrix} \begin{pmatrix} W_1 \xi \\ W_2 \xi \end{pmatrix}.$$

**Remark 2.3**. Lemma 2.5 is the alleged reciprocally convex inequality in the quaternion domain.

## 3. MAIN RESULTS

**Theorem 3.1**. Suppose Assumptions 1 and 2 hold. If there exist positive diagonal matrices $M_1, M_2, M_3 \in \mathbb{R}^{n \times n}$, positive definite matrices $P_\iota (\iota = 1,2,3), Q_\varsigma (\varsigma = 1,2,\ldots,6), R_1, R_2 \in \mathbb{Q}^{n \times n}$ and appropriate matrices $U, V, S_1, S_2 \in \mathbb{Q}^{n \times n}$ such that the following quaternion-valued LMIs hold:

$$\begin{pmatrix} R_1 & U \\ * & R_1 \end{pmatrix} > 0, \quad (3)$$

$$\begin{pmatrix} R_2 & V \\ * & R_2 \end{pmatrix} > 0, \quad (4)$$

$$\Omega = (\Omega_{ij})_{11 \times 11} < 0, \quad (5)$$

where

$$\Omega_{ji} = \Omega_{ij}^* (i \neq j), \Omega_{11} = -P_1 C - C P_1 + P_2 + \delta^2 P_3 + Q_1 + Q_3 + Q_5 + Q_6 - R_1 + \Gamma M_1 \Gamma, \Omega_{14}$$
$$= R_1 - U^*, \Omega_{16} = U^*, \Omega_{18} = P_1 A, \Omega_{1,10} = P_1 B, \Omega_{1,11} = C P_1 C, \Omega_{22}$$
$$= d_1^2 R_1 + d_2^2 R_2 - S_1 - S_1^*, \Omega_{23} = -S_1^* C - S_2, \Omega_{28} = S_1^* A, \Omega_{2,10} = S_1^* B, \Omega_{33}$$
$$= -P_2 - C S_2 - S_2^* C, \Omega_{38} = S_2^* A, \Omega_{3,10} = S_2^* B, \Omega_{44}$$
$$= -(1 - \mu_1)Q_1 - R_1 - R_1^* + U + U^* + \Gamma M_2 \Gamma, \Omega_{46} = R_1 - U^*, \Omega_{55}$$
$$= -(1 - \mu)Q_3 - R_2 - R_2^* + V + V^* + \Gamma M_3 \Gamma, \Omega_{56} = R_2^* - V, \Omega_{57}$$
$$= R_2 - V^*, \Omega_{66} = -Q_5 - R_1 - R_2, \Omega_{67} = V^*, \Omega_{77} = -Q_6 - R_2, \Omega_{88}$$
$$= Q_2 + Q_4 - M_1, \Omega_{8,11} = -A^* P_1 C, \Omega_{99} = -(1 - \mu_1)Q_2 - M_2, \Omega_{10,10}$$
$$= -(1 - \mu)Q_4 - M_3, \Omega_{10,11} = -B^* P_1 C, \Omega_{11,11} = -P_3.$$

Then the quaternion-valued neural networks is globally asymptotically stable.

**Proof.** Construct the following Lyapunov-Krasovskii functional

$$V(t) = \sum_{i=1}^{4} V_i(t) \quad (6)$$

Where

$$V_1(t) = \left(x(t) - C\int_{t-\delta}^{t} x(s)ds\right)^* P_1\left(x(t) - C\int_{t-\delta}^{t} x(s)ds\right), (7)$$

$$V_2(t) = \int_{t-\delta}^{t} x^*(s)P_2 x(s)ds + \delta \int_{-\delta}^{0}\int_{t+\theta}^{t} x^*(s)P_3 x(s)ds d\theta, (8)$$

$$V_3(t) = \int_{t-d_1(t)}^{t} (x^*(s)Q_1 x(s) + f^*(x(s))Q_2 f(x(s)))ds$$
$$+ \int_{t-d(t)}^{t} (x^*(s)Q_3 x(s) + f^*(x(s))Q_4 f(x(s)))ds + \int_{t-d_1}^{t} x^*(s)Q_5 x(s)ds$$
$$+ \int_{t-d}^{t} x^*(s)Q_6 x(s)ds, (9)$$

$$V_4(t) = d_1 \int_{-d_1}^{0}\int_{t+\theta}^{t} \dot{x}^*(s)R_1 \dot{x}(s)ds d\theta + d_2 \int_{-d}^{-d_1}\int_{t+\theta}^{t} \dot{x}^*(s)R_2 \dot{x}(s)ds d\theta. (10)$$

Then the derivatives of $V_i$ ($i = 1,2,3,4$) can be calculated and estimated straightforwardly:

$$\dot{V}_1(t) = -x^*(t)(P_1 C + CP_1)x(t) + x^*(t)P_1 Af(x(t)) + f^*(x(t))A^* P_1 x(t)$$
$$+ x^*(t)P_1 Bf\left(x(t-d(t))\right) + f^*\left(x(t-d(t))\right)B^* P_1 x(t)$$
$$+ \left(\int_{t-\delta}^{t} x(s)ds\right)^* CP_1 Cx(t) + x^*(t)CP_1 C\left(\int_{t-\delta}^{t} x(s)ds\right)$$
$$- \left(\int_{t-\delta}^{t} x(s)ds\right)^* CP_1 Af(x(t)) - f^*(x(t))A^* P_1 C\left(\int_{t-\delta}^{t} x(s)ds\right)$$
$$- \left(\int_{t-\delta}^{t} x(s)ds\right)^* CP_1 Bf\left(x(t-d(t))\right)$$
$$- f^*\left(x(t-d(t))\right)B^* P_1 C\left(\int_{t-\delta}^{t} x(s)ds\right), (11)$$

$$\dot{V}_2(t) \leq x^*(t)(P_2 + \delta^2 P_3)x(t) - x^*(t-\delta)P_2 x(t-\delta)$$
$$+ \left(\int_{t-\delta}^{t} x(s)ds\right)^* P_3\left(\int_{t-\delta}^{t} x(s)ds\right), (12)$$

$$\dot{V}_3(t) \leq x^*(t)(Q_1 + Q_3 + Q_5 + Q_6)x(t) - x^*(t - d_1)Q_5x(t - d_1) - x^*(t - d)Q_6x(t - d)$$
$$- (1 - \mu_1)x^*(t - d_1(t))Q_1x(t - d_1(t)) - (1 - \mu)x^*(t - d(t))Q_3x(t - d(t))$$
$$+ f^*(x(t))(Q_2 + Q_4)f(x(t))$$
$$- (1 - \mu_1)f^*(x(t - d_1(t)))Q_2f(x(t - d_1(t)))$$
$$- (1 - \mu)f^*(x(t - d(t)))Q_2f(x(t - d(t))), (13)$$

$$\dot{V}_4(t) = \dot{x}^*(t)(d_1^2 R_1 + d_2^2 R_2)\dot{x}(t) - d_1 \int_{t-d_1}^{t} \dot{x}^*(s) R_1 \dot{x}(s)ds - d_2 \int_{t-d}^{t-d_1} \dot{x}^*(s) R_2 \dot{x}(s)ds, (14)$$

where Lemma 2.4 has been applied in the estimate of $\dot{V}_2(t)$ in (12). Based on Lemma 2.5, we further estimate two integration terms in $\dot{V}_4(t)$ as

$$-d_1 \int_{t-d_1}^{t} \dot{x}^*(s) R_1 \dot{x}(s)ds = -d_1 \int_{t-d_1(t)}^{t} \dot{x}^*(s) R_1 \dot{x}(s)ds - d_1 \int_{t-d_1}^{t-d_1(t)} \dot{x}^*(s) R_1 \dot{x}(s)ds$$

$$\leq -\frac{d_1}{d_1(t)}\left(\int_{t-d_1(t)}^{t} \dot{x}(s)ds\right)^* R_1 \left(\int_{t-d_1(t)}^{t} \dot{x}(s)ds\right)$$

$$-\frac{d_1}{d_1 - d_1(t)}\left(\int_{t-d_1}^{t-d_1(t)} \dot{x}(s)ds\right)^* R_1 \left(\int_{t-d_1}^{t-d_1(t)} \dot{x}(s)ds\right)$$

$$\leq -\begin{pmatrix} x(t) - x(t - d_1(t)) \\ x(t - d_1(t)) - x(t - d_1) \end{pmatrix}^* \begin{pmatrix} R_1 & U \\ * & R_1 \end{pmatrix} \begin{pmatrix} x(t) - x(t - d_1(t)) \\ x(t - d_1(t)) - x(t - d_1) \end{pmatrix}$$

$$= -\begin{pmatrix} x(t) \\ x(t - d_1) \\ x(t - d_1(t)) \end{pmatrix}^* \begin{pmatrix} R_1 & -U^* & -R_1 + U^* \\ -U & R_1 & -R_1 + U \\ -R_1 + U & -R_1 + U^* & 2R_1 - U - U^* \end{pmatrix} \begin{pmatrix} x(t) \\ x(t - d_1) \\ x(t - d_1(t)) \end{pmatrix}. (15)$$

Analogously, one can obtain that

$$-d_2 \int_{t-d}^{t-d_1} \dot{x}^*(s) R_2 \dot{x}(s)ds = -d_2 \int_{t-d(t)}^{t-d_1} \dot{x}^*(s) R_2 \dot{x}(s)ds - d_2 \int_{t-d}^{t-d(t)} \dot{x}^*(s) R_2 \dot{x}(s)ds$$

$$\leq -\begin{pmatrix} x(t - d_1) \\ x(t - d) \\ x(t - d(t)) \end{pmatrix}^* \begin{pmatrix} R_2 & -V^* & -R_2 + V^* \\ -V & R_2 & -R_2 + V \\ -R_2 + V & -R_2 + V^* & 2R_2 - V - V^* \end{pmatrix} \begin{pmatrix} x(t - d_1) \\ x(t - d) \\ x(t - d(t)) \end{pmatrix}. (16)$$

In addition, it follows from Assumption 1 that

$$0 \leq x^*(t)\Gamma M_1 \Gamma x(t) - f^*(x(t))M_1 f(x(t)), (17)$$

$$0 \leq x^*(t - d_1(t))\Gamma M_2 \Gamma x(t - d_1(t)) - f^*(x(t - d_1(t)))M_1 f(x(t - d_1(t))), (18)$$

$$0 \leq x^*(t - d(t))\Gamma M_3 \Gamma x(t - d(t)) - f^*(x(t - d(t)))M_3 f(x(t - d(t))). (19)$$

By applying the free weighting matrix method, we gather from (3) that

$$0 = [S_1 \dot{x}(t) + S_2 x(t - \delta)]^* H + H^*[S_1 \dot{x}(t) + S_2 x(t - \delta)], (20)$$

Where

$$H = -\dot{x}(t) - Cx(t-\delta) + Af(x(t)) + Bf(x(t-d(t))). \quad (21)$$

By substituting (21) into (20), we proceed to obtain that

$$\begin{aligned}
0 = &-\dot{x}^*(t)(S_1 + S_1^*)\dot{x}(t) - \dot{x}^*(t)(S_1^*C + S_2)x(t-\delta) - x^*(t-\delta)(CS_1 + S_2^*)\dot{x}(t) \\
&+ \dot{x}^*(t)S_1^*Af(x(t)) + f^*(x(t))A^*S_1\dot{x}(t) + \dot{x}^*(t)S_1^*Bf(x(t-d(t))) \\
&+ f^*(x(t-d(t)))B^*S_1\dot{x}(t) - x^*(t-\delta)(CS_2 + S_2^*C)x(t-\delta) \\
&+ x^*(t-\delta)S_2^*Af(x(t)) + f^*(x(t))A^*S_2x(t-\delta) \\
&+ x^*(t-\delta)S_2^*Bf(x(t-d(t))) + f^*(x(t-d(t)))B^*S_2x(t-\delta). \quad (22)
\end{aligned}$$

Therefore, it follows from (11)-(19) and (22) that

$$\dot{V}(t) \leq \eta^*(t)\Omega\eta(t), \quad (23)$$

Where

$$\eta(t) = \big(x^*(t),\ \dot{x}^*(t),\ x^*(t-\delta),\ x^*(t-d_1(t)),\ x^*(t-d(t)),\ x^*(t-d_1),\ x^*(t-d),$$
$$f^*(x(t)),\ f^*(x(t-d_1(t))),\ f^*(x(t-d(t))),\ (\textstyle\int_{t-\delta}^{t} x(s)ds)^*\big)^*.$$

Then it follows from (5) and (23) that $\dot{V}(t) < 0$, which together with the radial unboundedness of $V(t)$ guarantee the global asymptotical stability of the QVNNs (2). The proof is completed.

**Remark 3.1**. Since the QVLMIs (3)-(5) cannot be straightforwardly resolved via the Matlab LMI toolbox, it is necessary to convert QVLMIs into CVLMIs to acquire a set of feasible solutions with the assistance of Lemmas 2.2 and 2.3. Based on Lemma 2.2, we first conduct plural decomposition on the quaternion matrix appeared in Theorem 3.1: $P_\iota = P_{\iota 1} + P_{\iota 2}j$ ($\iota = 1,2,3$), $Q_\varsigma = Q_{\varsigma 1} + Q_{\varsigma 2}j$ ($\varsigma = 1,2,\ldots,6$), $R_1 = R_{11} + R_{12}j, R_2 = R_{21} + R_{22}j, U = U_1 + U_2j, V = V_1 + V_2j, S_1 = S_{11} + S_{12}j, S_2 = S_{21} + S_{22}j$. Then the following corollary can be immediately obtained by resorting to Lemma 2.3.

**Corollary 3.1**. Suppose Assumptions 1 and 2 hold. The equilibrium of system (2) is globally asymptotically stable if there exist positive diagonal matrices, Hermitian matrices, skew symmetric matrices and matrices such that the following complex-valued LMIs hold:

$$\begin{pmatrix} P_{i1} & -P_{i2} \\ \bar{P}_{i2} & \bar{P}_{i1} \end{pmatrix} > 0, \quad \begin{pmatrix} Q_{j1} & -Q_{j2} \\ \bar{Q}_{j2} & \bar{Q}_{j1} \end{pmatrix} > 0, \quad (24)$$

$$\begin{pmatrix} X_1 & -X_2 \\ \bar{X}_2 & \bar{X}_1 \end{pmatrix} > 0, \quad \begin{pmatrix} Y_1 & -Y_2 \\ \bar{Y}_2 & \bar{Y}_1 \end{pmatrix} > 0, \quad \begin{pmatrix} \Omega_1 & -\Omega_2 \\ \bar{\Omega}_2 & \bar{\Omega}_1 \end{pmatrix} < 0, \quad (25)$$

Where

$$X_1 = \begin{pmatrix} R_{11} & U_1 \\ * & R_{11} \end{pmatrix}, \quad X_2 = \begin{pmatrix} R_{12} & U_2 \\ \blacksquare & R_{12} \end{pmatrix}, \quad Y_1 = \begin{pmatrix} R_{21} & V_1 \\ * & R_{21} \end{pmatrix}, \quad Y_2 = \begin{pmatrix} R_{22} & V_2 \\ \blacksquare & R_{22} \end{pmatrix}$$

and $\Omega_1 = (\Omega_{ij}^{(1)})_{11\times 11}, \Omega_2 = (\Omega_{ij}^{(2)})_{11\times 11}$ with $\Omega_{ij}^{(1)}, \Omega_{ij}^{(2)}$ omitted here due to the limited space.

## 4. A NUMERICAL EXAMPLE

In this section, an illustrative example is provided to validate the effectiveness of the theoretical results.
Consider the two-dimensional QVNNs (2) with parameters given as follows:

$$C = diag\{8,12\}, A = (a_{ij})_{2\times 2}, B = (b_{ij})_{2\times 2}, \delta = 0.5,$$

$$d_1(t) = 0.45\sin(t) + 0.25, d_2(t) = 0.15\cos(t) - 0.05, f_1(s) = f_2(s) = 0.2\tanh(s), s \in \mathbb{Q}$$

Where

$$\begin{aligned}
a_{11} &= 1.2 + 3.0i - 3.6j + 2.0k, & a_{12} &= 1.8 + 1.6i - 2.0j - 1.9k, \\
a_{21} &= 3.8 - 3.8i + 2.0j - 2.1k, & a_{22} &= 1.5 + 3.2i - 3.6j + 3.0k, \\
b_{11} &= 1.5 - 3.3i + 2.6j + 1.1k, & b_{12} &= 1.5 + 2.6i + 0.9j - 2.9k, \\
b_{21} &= 2.5 + 3.2i - 0.7j - 1.5k, & b_{22} &= 2.9 + 3.5i + 1.3j + 1.5k.
\end{aligned}$$

In accordance with Remark 2.1, the quaternion-valued matrices $A$ and $B$ can be uniquely expressed as $A = A_1 + A_2 j$ and $B = B_1 + B_2 j$ respectively, where

$$A_1 = \begin{pmatrix} 1.2 + 3.0i & 1.8 + 1.6i \\ 3.9 - 3.8i & 1.5 + 3.2i \end{pmatrix}, A_2 = \begin{pmatrix} -3.6 + 2.0i & -2.0 - 1.9i \\ 2.0 - 2.1i & -3.6 + 3.0i \end{pmatrix},$$

$$B_1 = \begin{pmatrix} 1.5 - 3.3i & 1.5 + 2.6i \\ 2.5 + 3.2i & 2.9 + 3.5i \end{pmatrix}, B_2 = \begin{pmatrix} 2.6 + 1.1i & 0.9 - 2.9i \\ -0.7 - 1.5i & 1.3 + 1.5i \end{pmatrix}.$$

In addition, it can be readily verified that Assumptions 1 and 2 are satisfied, and $d_1 = 0.7, d_2 = 0.1, d = 0.8, \mu_1 = 0.45, \mu_2 = 0.15, \mu = 0.6$. Therefore, a set of feasible solutions to CVLMIs (24)-(25) can be established via the Yalmip toolbox in Matlab (only partial matrices in the solutions are listed here due to the limited space):

$$M_1 = diag\{0.3635, 0.3544\}, M_2 = diag\{0.1306, 0.1306\}, M_3 = diag\{0.3254, 0.3341\},$$

$$P_{11} = \begin{pmatrix} 0.4610 + 0.0000i & 0.0012 - 0.0025i \\ 0.0012 + 0.0025i & 0.1332 + 0.0000i \end{pmatrix},$$

$$Q_{11} = \begin{pmatrix} -0.2628 + 0.0000i & 0.0001 - 0.0001i \\ 0.0001 + 0.0001i & -0.2561 + 0.0000i \end{pmatrix},$$

$$U_1 = \begin{pmatrix} -0.1404 + 0.0000i & -0.0001 + 0.0002i \\ -0.0001 - 0.0002i & -0.1423 - 0.0000i \end{pmatrix},$$

$$V_1 = \begin{pmatrix} -0.1698 + 0.0000i & -0.0000 + 0.0001i \\ -0.0000 - 0.0001i & -0.1686 - 0.0000i \end{pmatrix},$$

$$S_{11} = \begin{pmatrix} 0.3059 - 0.0000i & 0.0634 - 0.1002i \\ 0.0634 + 0.1001i & 0.2655 + 0.0000i \end{pmatrix},$$

$$S_{21} = \begin{pmatrix} 0.2465 + 0.0000i & 0.0766 - 0.1211i \\ 0.0511 + 0.0807i & 0.3208 + 0.0000i \end{pmatrix}.$$

Therefore, the equilibrium of the QVNNs (2) is globally asymptotically stable according to Corollary 1. Fig. 1 depicts the transient behavior of the neuron state in (2).

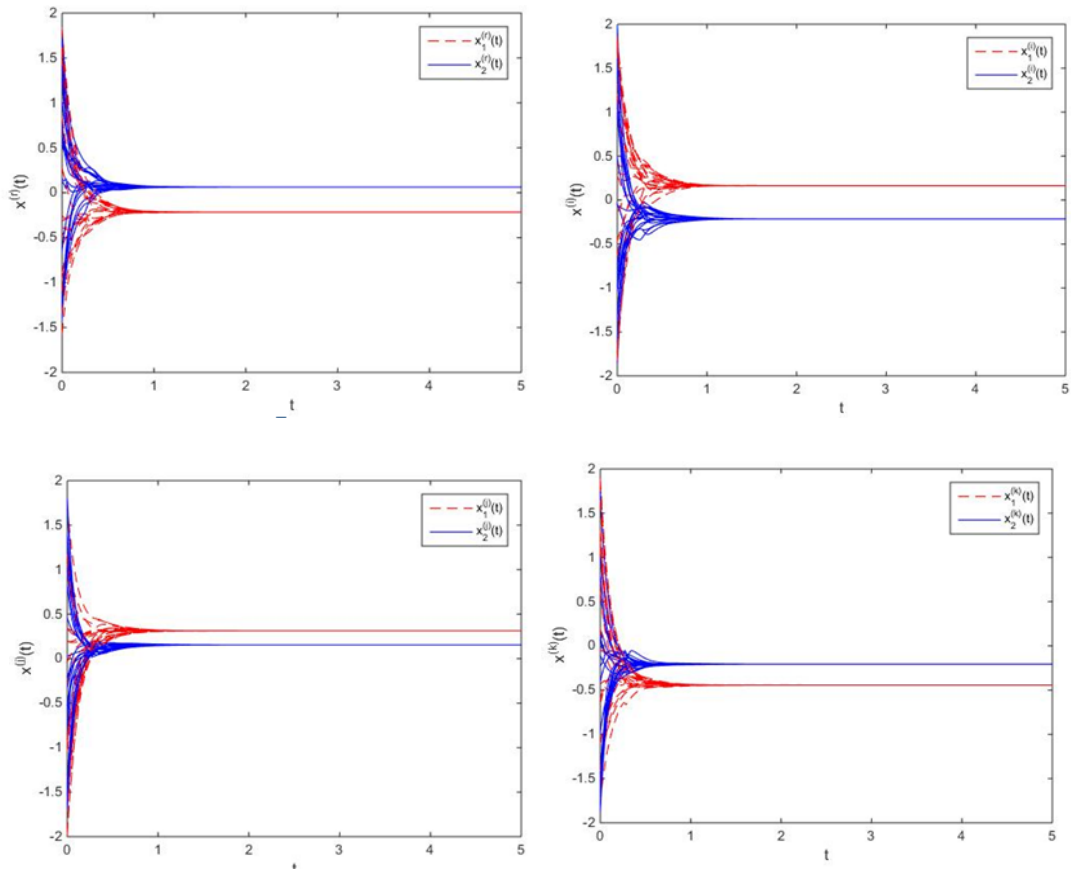

Figure 1. The transient behaviors of neuron states of QVNNs (2).

## 5. CONCLUSIONS

This paper is concerned with the stability analysis of quaternion-valued neural networks with both leakage delay and additive time-varying delays. Based on the Lyapunov functional method and inequality technique, some delay-dependent criteria are provided by fully considering the relationship between time-varying delays and upper bounds of delays. It is worth mentioning that the stability criteria are established in two forms. Finally, a numerical example is proposed to demonstrate the validity of theoretical results.

### ACKNOWLEDGEMENTS

The authors would like to thank the anonymous reviewers for their valuable comments and suggestions that help improve the quality of this manuscript.

### REFERENCES

[1] H. Bao, Ju H. Park, and J. Cao, (2016) "Exponential synchronization of coupled stochastic memristor-based neural networks with time-varying probabilistic delay coupling and impulsive delay", IEEE Transactions on Neural Networks, Vol. 27, No. 1, pp. 190–201.

[2] J. Cao and R. Li, (2017) "Fixed-time synchronization of delayed memristor-based recurrent neural networks", Science China: Information Sciences, Vol. 60, No. 3, pp. 108–122.

[3]  X. Li and J. Cao, (2017) "An impulsive delay inequality involving unbounded time-varying delay and applications",IEEE Transactions on Automatic Control, Vol. 62, No. 7, pp. 3618–3625.
[4]  J. Wang, H. Wu, and T. Huang, (2015) "Passivity-based synchronization of a class of complex dynamical networks with time-varying delay",Automatica, Vol. 56, pp. 105–112.
[5]  R. Yang, B. Wu, and Y. Liu, (2015) "A halanay-type inequality approach to the stability analysis of discrete time neural networks with delays", Applied Mathematics and Computation, Vol. 265, pp. 696–707.
[6]  X. Yang and J. Lu, (2016) "Finite-time synchronization of coupled networks with markovian topology and impulsive effects", IEEE Transactions on Automatic Control, Vol. 61, No. 8, pp. 2256–2261.
[7]  S. Jankowski, A. Lozowski, and J.M. Zurada, (1996) "Complex-valued multistate neural associative memory", IEEE Transactions on Neural Networks, Vol. 7, No. 6, pp. 1491–1496.
[8]  H. Bao, Ju H. Park, and J. Cao, (2016) "Synchronization of fractional-order complex-valued neural networks with time delay", Neural Networks, Vol. 81, pp. 16–28.
[9]  W. Gong, J. Liang, and J. Cao, (2015) "Matrix measure method for global exponential stability of complex-valued recurrent neural networks with time-varying delays",Neural Networks, Vol. 70, pp. 81–89.
[10] A. Hirose and S. Yoshida, (2012) "Generalization characteristics of complex-valued feedforward neural networks in relationto signal coherence", IEEE Transactions on Neural Networks and Learning Systems, Vol. 23, No. 4, pp. 541–551.
[11] R. Rakkiyappan, G. Velmurugan, X. Li, and D. O'Regan, (2016) "Global dissipativity of memristor-based complex-valued neural networks with time-varying delays", Neural Computing and Applications, Vol. 27, No. 3, pp. 629–649.
[12] N. Matsui, T. Isokawa, H. Kusamichi, F. Peper, and H. Nishimura, (2004) "Quaternion neural network with geometrical operators", Journal of Intelligent and Fuzzy Systems, Vol. 15, No. 3, pp. 149–164.
[13] B.C. Ujang, C.C. Took, and D.P. Mandic, (2011) "Quaternion-valued nonlinear adaptive filtering",IEEE Transactions onNeural Networks, Vol. 22, No. 8, pp. 1193–1206.
[14] Y. Liu, D. Zhang, J. Lu, and J. Cao, (2016) "Global μ-stability criteria for quaternion-valued neural networks with unbounded time-varying delays",Information Sciences, Vol. 360, pp. 273–288.
[15] H. Shu, Q. Song, Y. Liu, Z. Zhao, and F.E. Alsaadi, (2017) "Global μ-stability of quaternion-valued neural networks withnon-differentiable time-varying delays" ,Neurocomputing, Vol. 247, pp. 202–212.
[16] Y. Liu, D. Zhang, and J. Lu, (2016) "Global exponential stability for quaternion-valued recurrent neural networks with time-varying delays", Nonlinear Dynamics, Vol. 87, No. 1, pp. 1–13.
[17] X. Chen, Z. Li, Q. Song, J. Hu, and Y. Tan, (2017) "Robust stability analysis of quaternion-valued neural networks with time delays and parameter uncertainties", Neural Networks, Vol. 91, pp. 55–65.
[18] X. Chen, Q. Song, Z. Li, Z. Zhao, and Y. Liu, (2018) "Stability analysis of continuous-time and discrete-time quaternion-valued neural networks with linear threshold neurons", IEEE Transactions on Neural Networks, Vol. 29, No. 7, pp. 2769–2781.
[19] Z. Tu, J. Cao, A. Alsaedi, and T. Hayat, (2017) "Global dissipativity analysis for delayed quaternion-valued neural networks", Neural Networks, Vol. 89, pp. 97–104.
[20] Z. Yu, H. Gao, and S. Mou, (2008) "Asymptotic stability analysis of neural networks with successive time delay components", Neurocomputing, Vol. 71, No. 13-15, pp. 2848–2856.
[21] H. Shao and Q. Han, (2011) "New delay-dependent stability criteria for neural networks with two additive time-varying delay components", IEEE Transactions on Neural Networks, Vol. 22, No. 5, pp. 812–818.
[22] J. Tian and S. Zhong, (2012) "Improved delay-dependent stability criteria for neural networks with two additive time-varying delay components", Neurocomputing, Vol. 77, pp. 114–119.
[23] P.G. Park, J.W. Ko, and C. Jeong, (2011) "Reciprocally convex approach to stability of systems with time-varying delays", Automatica, Vol. 47, No. 1, pp. 235–238.
[24] J. Liang, K. Li, Q. Song, Z. Zhao, Y. Liu, and F.E. Alsaadi, (2018) "State estimation of complex-valued neural networks with two additive time-varying delays" ,Neurocomputing, Vol. 309, pp. 54–61.

**AUTHORS**

**Qun Huang** received the B.S. degree from Southeast University, Nanjing, China in 2014. He is currently pursuing the Ph.D. degree with the School of Mathematics, Southeast University,Nanjing, China. His current research interests include quaternion-valued neural networks and dynamic equations on time scales.

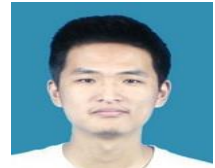

**Jinde Cao** (F'16) received the B.S. degree from Anhui Normal University, Wuhu, China, the M.S. degree from Yunnan University, Kunming, China, and the Ph.D. degree from Sichuan University, Chengdu, China, all in mathematics/applied mathematics, in 1986, 1989, and 1998, respectively. He is an Endowed Chair Professor, the Dean of the School of Mathematics, the Director of the Jiangsu Provincial Key Laboratory of Networked Collective Intelligence of China and the Director of the Research Center for Complex Systems and Network Sciences at Southeast University. Prof. Cao was a recipient of the National Innovation Award of China, Obada Prize and the Highly Cited Researcher Award in Engineering, Computer Science, and Mathematics by Thomson Reuters/Clarivate Analytics. He is a fellow of IEEE, a member of the Academy of Europe, a member of the European Academy of Sciences and Arts, a fellow of Pakistan Academy of Sciences, and an IASCYS academician.

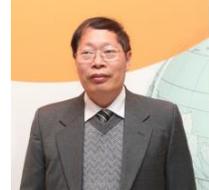